\documentclass{amsart}

\newtheorem{thm}{Theorem}[section]
\newtheorem{cor}[thm]{Corollary}
\newtheorem{lem}[thm]{Lemma}

\theoremstyle{definition}
\newtheorem{defn}[thm]{Definition}
\theoremstyle{remark}

\numberwithin{equation}{section}

\newcommand{\D}{\mathcal{D}}
\newcommand{\C}{\mathcal{C}}

\newcommand{\bfOmega}{\mathbf{\Omega}}
\newcommand{\bftau}{\mathbf{\tau}}
\newcommand{\DPRM}{\textsc{dprm }}
 
\newcommand{\abs}[1]{\left\vert#1\right\vert}
\newcommand{\set}[1]{\left\{#1\right\}}
\newcommand{\seq}[1]{\left\{#1\right\}}
 
\begin{document}

\title{On the existence of a new family of Diophantine equations for $\bfOmega$}
\author{Toby Ord}
\address{Department of Philosophy, The University of Melbourne, 3010, Australia}
\email{t.ord@pgrad.unimelb.edu.au}
\author{Tien D. Kieu}
\address{Centre for Atom Optics and Ultrafast Spectroscopy, Swinburne University of Technology,  Hawthorn 3122, Australia}
\email{kieu@swin.edu.au}
\maketitle



\begin{abstract}
We show how to determine the $k$-th bit of Chaitin's algorithmically random real number $\Omega$ by solving $k$ instances of the halting problem. From this we then reduce the problem of determining the $k$-th bit of $\Omega$ to determining whether a certain Diophantine equation with two parameters, $k$ and $N$, has solutions for an odd or an even number of values of $N$. We also demonstrate two further examples of $\Omega$ in number theory: an exponential Diophantine equation with a parameter $k$ which has an odd number of solutions iff the $k$-th bit of $\Omega$ is 1, and a polynomial of positive integer variables and a parameter $k$ that takes on an odd number of positive values iff the $k$-th bit of $\Omega$ is 1. 
\end{abstract}



\section{Introduction and Motivation}

In his 1975 paper \cite{Chaitin:1975}, Gregory Chaitin introduced the number $\Omega$. $\Omega$ is a real number between 0 and 1 which is the \emph{halting probability}: the probability that a randomly chosen program will halt. In order to make this precise, one must choose a particular (recursive) programming language in which the programs are self-delimiting binary strings. This allows us to choose a random program by flipping an unbiased coin to generate each bit, stopping when we reach a valid (self-delimited) program. This method makes longer programs increasingly unlikely and allows $\Omega$ to be a well defined probability. The exact value of $\Omega$ obviously depends upon which programming language we choose, and thus Chaitin~\cite{Chaitin:1987} specifies a particular language to make talk of $\Omega$ concrete. However, its important properties are the same regardless and the exact language chosen need not concern us.\footnote{The properties possessed by $\Omega$ are numerous and we can only discuss the most pertinent here. For further information on $\Omega$'s exotic attributes, see \cite{Chaitin:1987,Solovay:2000,Calude:2002,CaludeNies:1997}.}

It is well known that if we had access to $\Omega$ we could solve the halting problem~\cite{Chaitin:1987}. Furthermore, if we had some method of solving the halting problem then we could compute $\Omega$ to any desired accuracy. In this sense there is an equivalence between $\Omega$ and the halting problem: the ability to compute one gives us the ability to compute the other.

It is also well known that both the halting problem and $\Omega$ have analogues in number 
theory. The halting problem has been shown to be equivalent to the problem of determining
whether a given Diophantine equation has solutions, while $\Omega$ has been shown to be equivalent to determining whether a finite or infinite amount of Diophantine equations of a particular type have solutions.  However, unlike their counterparts in the theory of computability, these number theoretic versions have no direct connection between them. In this paper, we present a new representation of $\Omega$ within number theory that revolves around a matter of parity rather than finitude and makes clear the link between the number theoretic analogues of $\Omega$ and the halting problem.

In Sections~\ref{section:Randomness} and \ref{section:OmegaTau} we provide the formal definitions of $\Omega$ and a related real number $\tau$, presenting an efficient method for using the solutions to the halting problem to determine $\Omega$. In Sections~\ref{section:Diophantine} and \ref{section:OmegaFinitude}, we then introduce Diophantine equations and review the known method of using them to represent $\Omega$, providing proofs when they will facilitate the proofs of our new results. In Sections~\ref{section:OmegaParity} and \ref{section:Conclusions} we present and prove our results for a new representation of $\Omega$ and discuss their implications.


\section{Algorithmic Randomness and $\bfOmega$}
\label{section:Randomness}
\begin{defn}
$\Omega$ is formally defined by the following equation, where $p$ ranges over all programs in the language specified by Chaitin~\cite{Chaitin:1987} and $\abs{p}$ is the number of digits in the binary representation of $p$.
\begin{equation} \label{eq:Omega}
\Omega = \sum_{p \textrm{ halts}}2^{-\abs{p}}
\end{equation}
\end{defn}

The importance of $\Omega$ lies in its contrasting properties of being both an \emph{algorithmically random} and \emph{recursively enumerable} real number.

As a random real, $\Omega$ is uncomputable in a very strong sense. If we look at the binary expansion of $\Omega$ (the infinite sequence of 1's and 0's following the binary point), this sequence is highly incompressible. To see how this is so, let us compare it with the number $\tau$ from Copeland and Proudfoot~\cite{CopelandProudfoot:1999}, a non-recursive real that is not random.

\begin{defn}
$\tau$ is the real number between 0 and 1 whose $n$-th binary digit is 1 if the $n$-th program $p_n$ halts when given no input and 0 otherwise. In order to pick out a particular ordering of programs, let each program have a binary representation as above and let them be ordered in the usual lexical ordering for binary strings. We could also write:
\begin{equation}
\tau = \sum_{p_{n} \textrm{ halts}}2^{-n}
\end{equation}
\end{defn}

$\tau$ in an encoding of the answers to every instance of the halting problem in a single real number. While each bit of $\tau$ tells us whether a particular program halts when run with no input, we can easily use the bits of $\tau$ to answer the more common formulation of halting problem: `does program $p$ halt when run on input~$i$?'. Given $p$ and $i$, we can (recursively) construct a new program which takes no input and simulates the computation of $p$ on $i$. The bit of $\tau$ corresponding to this new program tells us whether or not $p$ halts on $i$.

The concept of compressibility for an infinite string can be understood through examining the minimum amount of advice needed to find out $n$ bits of that string~\cite{Chaitin:1987}. For recursive infinite strings (such as the binary digits of $\pi$), we can perform massive compression by asking for the bits of a particular program which will generate every bit of $\pi$, one by one. In this sense, there are only finitely many bits of information (those of that particular program) in the infinite bitstring that makes up $\pi$.

On the other hand, if we wish to determine the bits of $\tau$, there is no program that will do this, since the halting problem is undecidable. Therefore, determining $n$ bits of $\tau$ will require an amount of advice that increases with the value of $n$. Clearly we could find out $n$ bits of $\tau$ with $n$ bits of advice as we could be directly given those bits. However, Chaitin~\cite{Chaitin:1987} shows how we can do better than this, using only $\log_{2} n$ bits. This is because the $n$ bits of $\tau$ represent $n$ instances of the halting problem and to solve these, we only need to know how many of the $n$ programs halt. We can then simulate each program in parallel until this many have halted, being confident that every one that is still running will never halt. Therefore, while the (Turing) non-computability of $\tau$ makes it globally incompressible (computing infinitely many bits of $\tau$ requires infinitely many bits of advice), it is still locally compressible as $n$ bits of $\tau$ are computable from $\log_{2} n$ bits of advice.

In contrast, Chaitin~\cite{Chaitin:1987} defines a random real as one for which calculating $n$ bits of its binary expansion requires more than $n$ bits of advice. The reason that \emph{more} than $n$ bits are needed is that the advice, just like the program itself, must be self-delimiting. Chaitin~\cite{Chaitin:1975} has shown that $\Omega$ satisfies this condition and is thus a random real. In this way, the first $n$ bits of $\Omega$ contain $n$ bits of \emph{algorithmically incompressible} information.

In addition to recursive incompressibility, random reals are also characterised by \emph{recursive unpredictability}~\cite{Chaitin:1987}. Consider a `predictive' program that takes a finite initial segment of an infinite bitstring and returns a value indicating either `the next bit is 1', `the next bit is 0' or `no prediction'. If any such program is run on all finite prefixes of the binary expansion of a random real and makes an infinite amount of predictions, the limiting relative frequency of correct predictions approaches $\frac{1}{2}$. In other words when any program is used to predict infinitely many bits of a random real, such as $\Omega$, it does no better than random --- even with information about all the prior bits.
\\

In contrast to these results, $\Omega$ is a recursively enumerable (r.e.)~real. For a mathematical object to be r.e.~it need not be computable, but there must be a certain method of successively approximating it. For a set of positive integers to be r.e.~there must be a program that halts on $i$ if and only if $i$ is in the set. For an infinite binary sequence to be r.e.~there must be a program that halts on $i$ if and only if the $i$-th bit of the sequence is 1. For a real number $x$ to be r.e.~there must be a program that takes a positive integer $i$ and generates a rational approximation $x_i$ where $\seq{x_i}$ is an increasing sequence of rationals that converges to $x$. In the case of $\Omega$, we can construct such a sequence which we shall denote $\seq{\Omega_{i}}$, where
\begin{equation}
\Omega_i = \sum_{\stackrel{\scriptstyle \abs{p}\le i}{\scriptstyle p \textrm{ halts in } \le i \textrm{ steps}}} \!\!\!\!\!\!\!\!\!\!\!\!2^{-\abs{p}}
\label{OmegaN}
\end{equation} 

However, it is also important to ask whether the binary expansion of $\Omega$ is an r.e.~bitstring. For $\tau$, this is the case because there is a program that halts on $i$ if and only if the $i$-th bit is 1 --- indeed this is just a universal program, which simulates the $i$-th program (and thus halts iff the $i$-th program does). For $\Omega$ however, there is no such program. This is a rather surprising, but crucial, fact about $\Omega$.

We can see that this must be true given the incompressibility of $\Omega$. If there were a program to enumerate the bits of $\Omega$, we could use the trick described above for compressing $\tau$: we could determine $n$ bits of $\Omega$ by asking how many of these $n$ values make this program halt (and thereby require only $\log_2{n}$ bits). Therefore, as $\Omega$ is incompressible, we see that there is indeed no such program for $\Omega$ and that its sequence of bits is not r.e. Because all random reals share the property of being incompressible, this argument carries over and no random real can have an r.e.~sequence of bits.
\\

In the remainder of the paper, we shall use the fact that $\Omega$ is an r.e.~real to express it through Diophantine equations and thereby show how algorithmic randomness occurs even in the heart of number theory.


\section{Computing $\bfOmega$ from $\bftau$}
\label{section:OmegaTau}

Having seen that there are interesting connections between $\Omega$ and $\tau$, it is natural to ask whether we can compute $\Omega$ given access to $\tau$ (or in the language of recursion theory, whether $\Omega$ is \emph{Turing reducible} to $\tau$). One way to go about this is to construct a program $Q$ that takes two positive integers, $k$ and $N$, as input and halts iff $\frac{N}{2^{k}} < \Omega$. Such a program is quite easy to construct, it just needs to enumerate $\seq{\Omega_{i}}$ (as done in the previous section) and check at each stage whether $\frac{N}{2^{k}} < \Omega_{i}$, halting if this is true and continuing otherwise. Since $\seq{\Omega_{i}}$ approaches $\Omega$ from below, if $\frac{N}{2^{k}} < \Omega$ then there is some $i$ for which $\frac{N}{2^{k}} < \Omega_{i}$ and $Q$ halts as required. On the other hand, if $\frac{N}{2^{k}} \ge \Omega$ then there is no such $i$ and $Q$ will not halt. From $\tau$, we can determine whether or not $Q$ halts on a given input, and thus whether or not each $\frac{N}{2^{k}}$ is less than $\Omega$.

To determine the first $k$ bits of $\Omega$, we just need to determine the greatest value of $N$ for which $\frac{N}{2^{k}} < \Omega$. We will then know that for this $N$, $\frac{N}{2^{k}} < \Omega < \frac{N+1}{2^{k}}$ and thus that if $N$ is expressed as a $k$ digit binary number (with leading zeros if required), it will be the first $k$ digits of $\Omega$.

We can therefore use $\tau$ to determine the first $k$ bits of $\Omega$, by checking whether $Q$ halts when applied to $k$ and $N$ for each value of $N$ from 1 to $2^k-1$. This technique allows us to determine the first $k$ bits of $\Omega$ from $2^k-1$ carefully chosen bits of $\tau$. Indeed, given these bits, it is so simple to compute the first $k$ bits of $\Omega$, that we could do it by simply combining them with a truth table and $\Omega$ is not only Turing reducible to $\tau$, but also \emph{truth table reducible}.

However, if we are prepared to sacrifice this extreme algorithmic simplicity, we can get by with even fewer bits of $\tau$ by using a \emph{bisection search}. In the bits of $\tau$ that we need we know there is a lot of structure: if $\frac{N}{2^{k}} < \Omega$ then for all $i < N$, $\frac{i}{2^{k}} < \Omega$. Similarly, if $\frac{N}{2^{k}} \ge \Omega$ then for all $i > N$, $\frac{i}{2^{k}} \ge \Omega$. Thus, all we need to find is the greatest value of $N$ for which $\frac{N}{2^{k}} < \Omega$. This can be done by first trying $N=\frac{2^{k}}{2}$ and checking if it is less than or greater than $\Omega$. If it is greater, we then try the midpoint of 0 and $\frac{2^{k}}{2}$. If it is less, we try the midpoint of $\frac{2^{k}}{2}$ and $2^{k}$. We proceed in this manner until finding the value of $N$ for which $\frac{N}{2^{k}} < \Omega$ and $\frac{N+1}{2^{k}} \ge \Omega$.

Performing this type of bisection search on $n$ items requires only $\log_{2} (n+1)$ queries, rounding up to the nearest integer. As we do not need to test for $N=0$ or $N=2^{k}$, there are only $2^{k}-1$ remaining values to check and we can therefore get the first $k$ bits of $\Omega$ with only $k$ bits of $\tau$. From the incompressibility of $\Omega$, we can see that this must be optimal: no less than $k$ bits of any other real whatsoever can give us $k$ bits of $\Omega$.\footnote{It may seem like there is a contradiction here because we could use the compression trick for $\tau$ again, giving us a method for computing $k$ bits of $\Omega$ from only $\log_{2} k$ bits of advice. However, this turns out to be impossible because to use the bisection search, we need to find out the value of one bit of $\tau$ before we know which bit to ask for next. There is no time at which we know in advance a complete set of $k$ bits of $\tau$ that we can ask for in compressed form.}
\\

Having shown a method for determining the bits of $\Omega$ from an r.e.~sequence of bits ($\tau$), we will now examine Diophantine equations and look at how this method can be used to find $\Omega$ in number theory.


\section{Diophantine Equations and Hilbert's Tenth Problem}
\label{section:Diophantine}

\begin{defn}
A \emph{Diophantine equation} is an equation of the form
\begin{equation}
D(x_{1},\ldots,x_{m}) = 0
\end{equation}
where $D$ is a polynomial with integer coefficients. The variables are typically allowed to range over the integers or the non-negative integers, but here it will be convenient to restrict them to the positive integers.
\end{defn}

The study of Diophantine equations is a central area of number theory. It is well known that some Diophantine equations have solutions (such as $2x_{1} - 4 = 0$) and some do not (such as $2x_{1} - 3 = 0$). The task of devising an algorithm for determining whether or not an arbitrary Diophantine equation is solvable is known as Hilbert's tenth problem and has been a major area of research in 20th Century mathematics (see Matiyasevich~\cite{Matiyasevich:1993}).

\begin{defn}
A \emph{family of Diophantine equations} is a relation of the form
\begin{equation}
D(a_{1},\ldots,a_{n},x_{1},\ldots,x_{m}) = 0
\end{equation}
where we distinguish between the variables $a_{1},\ldots,a_{n}$ which are called parameters and $x_{1},\ldots,x_{m}$ which are called unknowns. Fixing values of the parameters, specifies one of the individual Diophantine equations that comprise the family. 
\end{defn}

\begin{defn}
A set of $n$-tuples, $\D$, is \emph{Diophantine} iff it can be represented by a family of Diophantine equations in the following way
\begin{equation}
\langle a_{1},\ldots,a_{n} \rangle  \in \D
\quad \Longleftrightarrow \quad
\exists x_{1}\ldots{}x_{m}{D(a_{1},\ldots,a_{n},x_{1},\ldots,x_{m})=0}
\end{equation}
\end{defn}

Many sets, such as the square numbers or the pairs $\langle a, b \rangle$ such that $a$ divides $b$ are Diophantine. The question of exactly which sets are Diophantine has received considerable attention and was solved through the work of Martin Davis, Hilary Putnam, Julia Robinson and Yuri Matiyasevich~\cite{Matiyasevich:1993}.

\begin{thm}[\DPRM Theorem]
Every r.e.~set $\D$ of $n$-tuples of positive integers is Diophantine, that is:
\begin{equation}
\langle a_{1},\ldots,a_{n} \rangle  \in \D
\quad \Longleftrightarrow \quad
\exists x_{1}\ldots{}x_{m}{D(a_{1},\ldots,a_{n},x_{1},\ldots,x_{m})=0}
\end{equation}
\end{thm}
\begin{proof}
See~\cite{Matiyasevich:1993}.
\end{proof}

The \DPRM Theorem provides an excellent tool for showing the existence of Diophantine equations with certain properties. All that is needed is to come up with a set that has the desired property and prove that it is r.e. For example, the following is immediately evident:

\begin{cor}[Family of Diophantine equations for $\bftau$]
There is a family of Diophantine equations:
\begin{equation}
T(n,x_{1},\ldots,x_{m})=0
\end{equation}
which has solutions for a given value of $n$ iff the $n$-th bit of $\tau$ is 1 or, equivalently, iff the $n$-th program halts.
\end{cor}
\begin{proof}
The set of positive integers $n$ such that the $n$-th bit of $\tau$ is 1 is clearly an r.e.~set. Thus, by the \DPRM Theorem, there is a family of Diophantine equations with a parameter, $n$, which has a solution for a given value of $n$ iff the $n$-th bit of $\tau$ is 1.
\end{proof}

From this it is clear that Hilbert's tenth problem must be recursively undecidable. No program could decide whether a given Diophantine equation has a solution because this would allow a program to compute the bits of $\tau$ --- a task that is known to be non-recursive.%
\footnote{Indeed, it was long known that the recursive undecidability of Hilbert's Tenth Problem would follow immediately from the \DPRM Theorem and this was the main motivation for its proof --- the Diophantine representations for all other r.e.~sets being a nice corollary.}
We therefore have undecidability in number theory, with direct analogues of the halting problem and $\tau$. In what follows we show how we can also find algorithmic randomness and, in doing so, add $\Omega$ to this list.


\section{Chaitin's Expression of $\bfOmega$ Through Diophantine Equations}
\label{section:OmegaFinitude}

From the \DPRM Theorem, we can see that $\Omega$ cannot be directly represented by a Diophantine set (one that has $k$ as a member iff the $k$-th bit of $\Omega$ is 1) because the bits of $\Omega$ are not r.e. However, taking a slightly different approach, Chaitin~\cite{Chaitin:1987} shows that the bits of $\Omega$ can be found in another property of Diophantine equations.

\begin{thm}[Family of Diophantine Equations for $\Omega$ via Finitude]
\label{theorem:FinitudeFamily}
There is a family of Diophantine equations 
\begin{equation}
\chi_{1}(k,N,x_{1},\ldots,x_{m})=0 
\end{equation}
such that for a given $k$ there are infinitely many values of $N$ for which it has a solution if the $k$-th bit of $\Omega$ is 1 and finitely many values of $N$ for which it has a solution if the $k$-th bit of $\Omega$ is 0.
\end{thm}
\begin{proof}
There is a program $P$ that takes two inputs, $k$ and $N$, generating $\Omega_{N}$ and returning the value of its $k$-th bit. As $N$ increases, the $k$-th bit of $\Omega_{N}$ may change between 0 and 1 many times, but since $\seq{\Omega_{N}}$ is monotonically increasing, each bit can only change finitely many times (at most $2^k-1$). Therefore, the $k$-th bit of $\Omega_{N}$ will eventually settle on either a 0 or a 1 and because $\seq{\Omega_{N}}$ converges to $\Omega$, this final value must be the same as the $k$-th bit of $\Omega$. Thus, the $k$-th bit of $\Omega_{N}$ will only differ from the $k$-th bit of $\Omega$ for finitely many values of $N$. 

Let $\C_{1}$ be the set of pairs $\langle k, N \rangle$ such that $P$ returns 1 when applied to $k$ and $N$. If the $k$-th bit of $\Omega$ is 0, $P$ will return 1 (the correct value for $\Omega_N$ but incorrect for $\Omega$) for only finitely many values of $N$.  On the other hand, if the $k$-th bit of $\Omega$ is 1, $P$ will return 1 for infinitely many values of $N$. In other words, for a given $k$, there are infinitely many pairs $\langle k,N \rangle \in \C_{1}$ if the $k$-th bit of $\Omega$ is 1 and finitely many if it is 0. Since $\C_{1}$ is recursive, it is r.e.~and thus via the \DPRM Theorem there is a Diophantine equation, $\chi_{1}(k,N,x_{1},\ldots,x_{m})=0$, that has a solution iff $\langle k,N\rangle \in \C_{1}$. Thus, for a given value of $k$, this equation has solutions for infinitely many values of $N$ iff the $k$-th bit of $\Omega$ is 1.
\end{proof}

Chaitin takes this result further and provides a simpler example of $\Omega$ occurring in number theory. In doing so he moves to \emph{exponential Diophantine equations}. Where the polynomial in a Diophantine equation consists of variables and integer constants composed together with addition and multiplication, exponential Diophantine equations merely add exponentiation to this list, with the proviso that no negative constants can appear in the exponents. In proving the \DPRM Theorem, Matiyasevich~\cite{Matiyasevich:1993} showed that all sets which are exponential Diophantine are Diophantine as well. However, Chaitin's move to exponential Diophantine equations allows the use of another key result due to Matiyasevich~\cite{Matiyasevich:1993} which concerns not just the existence of solutions to an equation, but the quantity of solutions as well.

\begin{thm}[Existence of Singlefold Exponential Diophantine Equations]
\label{singlefold}
For any Diophantine set, $\D$, there is a family of exponential Diophantine equations
\begin{equation}
E(a_{1},\ldots,a_{n},x_{1},\ldots,x_{m}) = 0
\end{equation}
that is a \emph{singlefold} representation of $\D$: it has exactly one solution for $\langle a_{1},\ldots,a_{n} \rangle \in \D$ and no solutions otherwise.
\end{thm}
\begin{proof}
See~\cite{Matiyasevich:1993}.
\end{proof}

\begin{thm}[Family of Exponential Diophantine Equations for $\bfOmega$ via Finitude]
\label{theorem:FinitudeExpFamily}
There is a family of exponential Diophantine equations 
\begin{equation}
\chi^e_{1}(k,x_{0},x_{1},\ldots,x_{m}) = 0
\end{equation}
with infinitely many solutions for $k$ if the $k$-th bit of $\Omega$ is 1 and finitely many if the $k$-th bit of $\Omega$ is 0.
\end{thm}

\begin{proof}
Let $\chi^e_{1}(k,N,x_{1},\ldots,x_{m}) = 0$ be the family of singlefold exponential Diophantine equations for $\chi_{1}$ of Theorem~\ref{theorem:FinitudeFamily}.
The existence of $\chi^e_1$ is ensured by Theorem~\ref{singlefold}.  
Given a value of $k$, there is a single solution for each of infinitely many values of $N$ iff the $k$-th bit of $\Omega$ is 1. If we treat $N$ as another unknown instead of a parameter, we get $\chi^e_{1}(k,x_{0},x_{1},\ldots,x_{m}) = 0$ which is a family of exponential Diophantine equations with infinitely many solutions for $k$ iff the $k$-th bit of $\Omega$~is~1.
\end{proof}

Thus, while the question of whether a Diophantine equation has solutions is undecidable  in general, the question of whether an exponential Diophantine equation has infinitely many solutions is much worse. The task is no longer r.e.~(since solving it gives the bits of $\Omega$) and as a single parameter is varied, the results fluctuate in an algorithmically random manner: there is absolutely no recursive pattern to be found.

Chaitin~\cite{Chaitin:1987} went even further than we have here by actually constructing an exponential Diophantine equation $\chi^e_{1}(k,x_{0},x_{1},\ldots,x_{m}) = 0$. This equation was automatically generated from a complex register machine program and is very large, with approximately 17,000 unknowns. While it is has been shown that this can be reduced to just three \cite{Matiyasevich:1993}, doing so would be a very challenging task.


\section{A New Expression of $\bfOmega$ Through Diophantine Equations}
\label{section:OmegaParity}

Using the method of computing $\Omega$ from $\tau$ discussed in Section \ref{section:OmegaTau} we can now present our main results. While the bits of $\Omega$ occur in Theorems \ref{theorem:FinitudeFamily} and \ref{theorem:FinitudeExpFamily} through the distinction between the finite and the infinite, the presentation below remains within the realms of the finite, with the bits of $\Omega$ occurring in the changes of parity.

\begin{thm}[Family of Diophantine Equations for $\bfOmega$ via Parity]
\label{theorem:ParityFamily}
There is a family of Diophantine equations
\begin{equation}
\chi_{2}(k,N,x_{1},\ldots,x_{m})=0 
\end{equation}
such that for a given $k$: 
\begin{itemize}
\item there are less than $2^{k}$ values of $N$, taken from $\set{1,\ldots,2^k-1}$, for which there is a solution,
\item and there are an odd number of values of $N$ for which there is a solution iff the $k$-th bit of $\Omega$ is 1.
\end{itemize}
\end{thm}

\begin{proof}
Let $\C_{2}$ be the set of pairs of positive integers $\langle k,N \rangle$ such that $\frac{N}{2^{k}} < \Omega$. From Section \ref{section:OmegaTau}, we see that $\C_{2}$ is r.e.~and so, by the \DPRM Theorem, there exists a Diophantine equation, $\chi_{2}(k,N,x_{1},\ldots,x_{m})=0$, which has a solution iff $\langle k,N\rangle \in \C_2$. Since $0 < \Omega < 1$, $\langle k,N \rangle \in \C_{2}$ implies that $N \in \set{1,\ldots,2^k-1}$, and so these are the only values of $N$ for which $\chi_{2}(k,N,x_{1},\ldots,x_{m})=0$ can have solutions.

If the greatest value of $N$ for which $\frac{N}{2^{k}} < \Omega$ is odd, then the $k$-th bit of $\Omega$ must be 1 and similarly, if the greatest value of $N$ for which $\frac{N}{2^{k}} < \Omega$ is even, then the $k$-th bit of $\Omega$ must be 0. As the greatest value of $N$ for which $\frac{N}{2^{k}} < \Omega$ is equal to the number of values of $N$ for which $\frac{N}{2^{k}} < \Omega$, it follows that for a given $k$, there are an odd number of values of $N$ for which $\chi_{2}(k,N,x_{1},\ldots,x_{m})=0$ has solutions iff the $k$-th bit of $\Omega$ is 1.
\end{proof}

Unlike the representation of $\Omega$ through finitude, this representation would allow $\Omega$ to be directly computed by determining whether certain Diophantine equations have solutions (i.e.~by solving several instances of Hilbert's tenth problem). This mirrors the relationship between $\Omega$ and the halting problem in computability theory. Even the efficiency of the computation is preserved: the bisection approach could also be used here to determine the $k$-th bit of $\Omega$ from the solutions to just $k$ instances of Hilbert's tenth problem.

\begin{thm}[Family of Exponential Diophantine Equations for $\bfOmega$ via Parity]
\label{theorem:ParityExpFamily}
There is a family of exponential Diophantine equations
\begin{equation}
\chi^e_{2}(k,x_{0},x_{1},\ldots,x_{m}) = 0
\end{equation}
which for a given $k$, 
\begin{itemize}
\item has less than $2^{k}$ solutions, where $x_0$ takes distinct values from
the set $\set{1,\ldots,2^k-1}$,
\item and has an odd number of solutions iff the $k$-th bit of $\Omega$ is 1.
\end{itemize}
\end{thm}

\begin{proof}
Let $\chi^e_{2}(k,N,x_{1},\ldots,x_{m}) = 0$ be the singlefold family of exponential Diophantine equations for $\C_{2}$. Given a value of $k$, there is a single solution for each of an odd number of values of $N$ iff the $k$-th bit of $\Omega$ is 1. If we treat $N$ as another unknown instead of as a parameter, we get $\chi^e_{2}(k,x_{0},x_{1},\ldots,x_{m}) = 0$.
For a given value of $k$, this has less than $2^{k}$ solutions and in these solutions, $x_{0}$ takes distinct values from $\set{1,\ldots,2^k-1}$. Furthermore, the number of solutions is odd iff the $k$-th bit of $\Omega$ is 1.
\end{proof}

In addition to the above results, we can use both the finitude and parity methods to produce a new source of algorithmic randomness in number theory. Using the following lemma~\cite{Matiyasevich:1993}, we can design polynomial expressions with variables and a parameter $k$ ranging over the positive integers that exhibit randomness in the number of positive values they assume as $k$ is varied.

\begin{lem}
\label{lemma:PolyLem}
Given a family of Diophantine equations with two parameters
\begin{equation}
D(k,N,x_{1},\ldots,x_{m}) = 0
\end{equation}
we can construct a polynomial with integer coefficients, one parameter and one additional variable $x_{0}$ as follows
\begin{equation}
P(k,x_{0},x_{1},\ldots,x_{m})
\equiv
x_{0}\left(1-(D(k,x_{0},x_{1},\ldots,x_{m}))^{2}\right).
\end{equation}
If we restrict the variables to values in the positive integers, then the set of positive values that this polynomial assumes for a given $k$ is exactly the set of all $N$ such that $D(k,N,x_{1},\ldots,x_{m}) = 0$ has solutions.
\end{lem}

\begin{proof}
For all $x_{0},x_{1},\ldots,x_{m}$ such that $D(k,x_{0},x_{1},\ldots,x_{m}) = 0$, the polynomial $P(k,x_{0},x_{1},\ldots,x_{m})$ will assume the value $x_{0}$. For all $x_{0},x_{1},\ldots,x_{m}$ such that $D(k,x_{0},x_{1},\ldots,x_{m}) \neq 0$, $P(k,x_{0},x_{1},\ldots,x_{m})$ will assume a non-positive value. Thus, $P(k,x_{0},x_{1},\ldots,x_{m})$ will assume all values of $N$ such that $D(k,N,x_{1},\ldots,x_{m}) = 0$ has solutions and no other positive values.
\end{proof}

\begin{thm}[Polynomial Expression for $\bfOmega$ via Finitude]
\label{theorem:FinitudePoly}
\begin{equation}
P_{1}(k,x_{0},x_{1},\ldots,x_{m})
\equiv
x_{0}\left(1-(\chi_{1}(k,x_{0},x_{1},\ldots,x_{m}))^{2}\right)
\end{equation}
is a polynomial with integer coefficients and a parameter $k$. If we restrict the variables to values in the positive integers, then this polynomial assumes infinitely many positive values iff the $k$-th bit of $\Omega$ is 1.
\end{thm}

\begin{proof}
By Lemma \ref{lemma:PolyLem}, the set of positive values that $P_{1}$ assumes for a given $k$ is precisely the set of $N$ such that $\chi_{1}(k,N,x_{1},\ldots,x_{m})=0$ has solutions. By Theorem \ref{theorem:FinitudeFamily}, this set has infinitely many members iff the $k$-th bit of $\Omega$ is 1, so $P_{1}$ must assume a countable infinity of positive values iff the $k$-th bit of $\Omega$ is 1.
\end{proof}

\begin{thm}[Polynomial Expression for $\bfOmega$ via Parity]\label{theorem:ParityPoly}
\begin{equation}
P_2(k,x_{0},x_{1},\ldots,x_{m})
\equiv
x_{0}\left(1-(\chi_2(k,x_{0},x_{1},\ldots,x_{m}))^{2}\right)
\end{equation}
is a polynomial with integer coefficients and a parameter $k$. If we restrict the variables to values in the positive integers, then 
\begin{itemize}
\item this polynomial assumes less than $2^k$ values in the positive integers,
taken from $\set{1,\ldots,2^k-1}$, 
\item and the number of positive values that it takes is odd iff the $k$-th bit of $\Omega$ is 1.
\end{itemize}
\end{thm}

\begin{proof}
From Lemma \ref{lemma:PolyLem}, the set of positive integer values that $P_{2}(k,x_{0},x_{1},\ldots,x_{m})$ assumes for a given $k$ is precisely the set of $N$ such that $\langle k,N \rangle \in \C_{2}$. 
By Theorem \ref{theorem:ParityFamily}, this set only includes values from $\set{1,\ldots,2^k-1}$. Theorem \ref{theorem:ParityFamily} also tells us that there are an odd number of $N$ such that $\langle k,N \rangle \in \C_{2}$ iff the $k$-th bit of $\Omega$ is 1, so $P_{2}(k,x_{0},x_{1},\ldots,x_{m})$ must assume an odd number of positive values iff the $k$-th bit of $\Omega$ is 1. 
\end{proof}

While Theorems \ref{theorem:ParityFamily}, \ref{theorem:ParityExpFamily} and \ref{theorem:ParityPoly} are presented in terms of producing the $k$-th bit of $\Omega$, the ideas behind them have a somewhat more general application. In each theorem, a set of values is discussed and its size is considered. For Theorem \ref{theorem:ParityFamily}, this is the number of values of $N$ for which $\chi_2 = 0$ has solutions, for Theorem \ref{theorem:ParityExpFamily} it is the number of solutions for $\chi^e_2 = 0$ and for Theorem \ref{theorem:ParityPoly} it is the number of positive values taken by $P_2$. For a given value of $k$, each of these quantities will equal the same number which we shall call $q_k$.

Previously, we just looked at the parity of $q_k$ and used this single bit of information to determine a single bit of $\Omega$. Since $q_k$ may take on $2^k$ different values, it contains $k-1$ additional bits of information and these can directly provide all prior bits of $\Omega$. When $q_k$ is expressed in binary (with enough leading zeros to give $k$ digits) it comprises the first $k$ bits of $\Omega$. For example, if for $k=5$, there are six solutions to $\chi^e_2 = 0$, then the first five bits of $\Omega$ are 00110.

In this way, we can see that the first $k-1$ bits of $q_k$ are just the $k-1$ bits of $q_{k-1}$. So, given $q_{k-1}$, there are only two possibilities for $q_k$: either $q_k = 2q_{k-1}$ or $q_k = 2q_{k-1}+1$. Thus, while the sequence $\seq{q_k \mod 2}$ is algorithmically random, the sequence $\seq{q_k}$ is highly structured.

We can also extend our results to the expression of $\Omega$ in different bases. If we consider the representation of $\Omega$ in base $b$, we can generate an equation for finding its $k$-th digit by replacing all references to $2^{k}$ with $b^{k}$. Instead of checking the parity of the appropriate quantity, the $k$-th digit is simply the remainder when this quantity is divided by $b$. Similarly, if we represent the quantity in base $b$ with enough leading zeros, it forms the first $k$ digits in a base $b$ representation of $\Omega$ and $q_k$ must be one of $\set{bq_{k-1}, bq_{k-1}+1, \ldots, bq_{k-1}+(b-1)}$.


\section{Concluding Remarks}
\label{section:Conclusions}
In this paper we report on our findings on the existence of a new family of Diophantine 
equations for determining the digits of $\Omega$.  Our equations differ from the previous 
families of Diophantine equations in several important features.  Previously, a particular bit of $\Omega$ was represented by whether some Diophantine equations have solutions for a finite or infinite number of values of a given parameter or, in the exponential Diophantine case, whether the number of solutions is finite. In contrast, our formulation always has a finite number of parameter values for which solutions could occur and a finite number of solutions in the case of the exponential Diophantine equations. The value of the corresponding bit of $\Omega$ depends only on whether this number is even or odd. Thus, in the relatively mundane switching between an odd or an even number of solutions as a parameter is varied, the full subtlety of algorithmic randomness is felt: the incompressibility and the unpredictability.

These new families of Diophantine equations for $\Omega$ also provide a symmetry to the relationship between randomness and undecidability in the fields of computability and number theory. Just as we can directly compute $\Omega$ (in its guise as the halting probability) from solutions to instances of the halting problem, so can we compute $\Omega$ (in its guise as a property of a family of Diophantine equations) from solutions to instances of Hilbert's tenth problem. Furthermore, the translation into the domain of number theory preserves the efficiency of the computation, producing $n$ bits of $\Omega$ from either the solutions to $n$ instances of the halting problem or the solutions to $n$ instances of Hilbert's tenth problem.

Finally, we provide a further example of how algorithmic randomness occurs in number theory: a polynomial with positive integer variables and a parameter $k$ that takes on an odd number of positive values iff the $k$-th bit of $\Omega$ is 1. Along with the other results of this paper, this hints at the variety of ways in which this randomness can occur. While the polynomials and exponentials that give rise to algorithmic randomness are not of the type that are likely to occur in classical research in number theory, their importance lies in showing that in some places there are facts with no recursive pattern at all. This is not to say that these facts are \emph{completely} patternless --- on the contrary, we have shown intricate dependencies between the expressions of $\Omega$ and those of the halting problem --- but in the face of so subtle a pattern, computer programs can perform no better than coin tossing. Without some means to transcend the fundamental limits of our current computers, whether through some kind of mathematical insight or radical new technology\footnote{Examples and discussions of such non-classical computation have been presented by both authors~\cite{Kieu:2002a,Ord:2002} and a very thorough survey can be found in Copeland~\cite{Copeland:2002}.}, they will remain completely beyond our grasp.

\section*{Acknowledgements}
We would like to thank Gregory Chaitin for helpful discussion.  TDK wishes to acknowledge 
the continuing support of Peter Hannaford.

\bibliographystyle{fundam}
\bibliography{diophantine,ait,hypercomputation}

\end{document}